\documentclass[11pt]{article}
\usepackage[usenames]{color}
\usepackage{amssymb,amsmath,psfig,epsfig,latexsym,graphics,here}
\usepackage[colorlinks=true,
linkcolor=webgreen,
filecolor=webbrown,
citecolor=webgreen]{hyperref}

\definecolor{webgreen}{rgb}{0,.5,0}
\definecolor{webbrown}{rgb}{.6,0,0}

\newtheorem{thm}{Theorem}

\newtheorem{cor}[thm]{Corollary}

\newenvironment{proof}{
        \noindent
        {\bf Proof:}
}{
        \hfill$\blacksquare$
}

\newcommand{\eeq}{\end{equation}}
\newcommand{\beql}[1]{\begin{equation}\label{#1}}
\newcommand{\eqn}[1]{(\ref{#1})}

\begin{document}
\begin{center}

{\Large{\bf On Non-Squashing Partitions }} \\
\vspace{1\baselineskip}
\vspace{1\baselineskip}
{\em N. J. A. Sloane}\footnote{Corresponding author. 
Postal address: AT\&T Shannon Labs, 180 Park Ave., Room C233,
Florham Park, NJ 07932--0971, USA;
phone: 973 360 8415; fax: 973 360 8718 }\\
AT\&T Shannon Labs \\
Information Sciences Research Center \\
Florham Park, NJ 07932--0971, USA \\
(Email address: \href{mailto:njas@research.att.com}{{\tt njas@research.att.com}})

\vspace{.5\baselineskip}
{\em James A. Sellers } \\
Department of Mathematics \\
Penn State University \\
University Park, PA 16802--6401, USA \\
(Email address: \href{mailto:sellersj@math.psu.edu }
{{\tt sellersj@math.psu.edu}})

\vspace{2\baselineskip}
{\bf Abstract} \\
\vspace{.4\baselineskip}
\end{center}

A partition $n = p_1 + p_2 + \cdots + p_k$ with
$1 \le p_1 \le p_2 \le \cdots \le p_k$
is called {\em non-squashing} if
$p_1 + \cdots + p_j \le p_{j+1}$ for 
$1 \le j \le k-1$.
Hirschhorn and Sellers showed that the number of
non-squashing partitions of $n$ is equal to the number
of binary partitions of $n$.
Here we exhibit an explicit bijection between
the two families, and
determine the number of non-squashing partitions
with distinct parts, with a specified
number of parts, or with a specified maximal part.
We use the results to solve a certain box-stacking problem.

\begin{center}{Keywords: partitions, non-squashing partitions, binary partitions, m-ary partitions, stacking boxes }\end{center}
\begin{center}{AMS 2000 Classification: Primary 11P81, 05A15}\end{center}
\vspace{.4\baselineskip}

\section{Introduction}

A correspondent, Claudio Buffara, recently asked for the solution
to the following problem, originally proposed by Telmo Luis Correia Jr.
We are given $n$ boxes, labeled $1, 2, \ldots, n$. For $i = 1, \ldots, n$,
box $i$ weighs $i$ grams and can support a total weight of $i$ grams.
What is $f(n)$, the number of different ways to
build a single stack of boxes in which no box will
be squashed by the weight of the boxes above it?
For example, $f(4) = 14$, since we
can form the following stacks:
$$
\emptyset \, ,   \,
1 \, ,  \,
2 \, ,  \,
3 \, ,  \,
4 \, ,  \,
\stackrel{\textstyle{1}}{2} \, ,  \,
\stackrel{\textstyle{1}}{3} \, ,  \,
\stackrel{\textstyle{1}}{4} \, ,  \,
\stackrel{\textstyle{2}}{3} \, ,  \,
\stackrel{\textstyle{2}}{4} \, ,  \,
\stackrel{\textstyle{3}}{4} \, ,  \,
\stackrel{\textstyle{1}}{\stackrel{\textstyle{2}}{3}} \, ,  \,
\stackrel{\textstyle{1}}{\stackrel{\textstyle{2}}{4}} \, ,  \,
\stackrel{\textstyle{1}}{\stackrel{\textstyle{3}}{4}} ~.
$$
The other two possible stacks:
$$
\stackrel{\textstyle{2}}{\stackrel{\textstyle{3}}{4}} \, ,~
\stackrel{\textstyle{1}}{\stackrel{\textstyle{2}}{\stackrel{\textstyle{3}}{4}}}  ~,
$$
are excluded, since $2+3 > 4$ and the box labeled $4$
would collapse in both cases.

To make this more precise, let us say that
a partition of a natural number $m$ into $k$ parts
is  {\em non-squashing} if when the parts are arranged in 
nondecreasing order, say
\beql{EqPart}
m = p_1 + p_2 + \cdots + p_k \mbox{~with~}
1 \le p_1 \le p_2 \le \cdots \le p_k ~,
\eeq
we have
\beql{EqNS}
p_1 + \cdots + p_j \le p_{j+1} \mbox{~for~} 
1 \le j \le k-1 ~.
\eeq
If the boxes in a stack are labeled (from
the top) $p_1, p_2, \ldots, p_k$, the stack will not
collapse if and only if the partition is
non-squashing. In the problem as stated, the boxes must
also have distinct labels and their sum cannot exceed
$\binom{n+1}{2}$. Therefore $f(n)$ is equal
to the total number of partitions of numbers from $0$ to
$\binom{n+1}{2}$ which are (i) non-squashing,
(ii) have distinct parts, and (iii) involve
no part greater than $n$.
We will give the solution in Sections \ref{BOXES} and \ref{BOXES2}.
In Sections \ref{Sec:NS} and \ref{NSPD}
we study the numbers of non-squashing partitions
and non-squashing partitions with distinct parts.
Sections \ref{NSNOP}, \ref{NSkDP} and \ref{NSLP}
deal with non-squashing partitions
with a given number of parts,
a given number of distinct parts,
and a specified largest part, respectively.
Some of these results are used in the final two sections,
others are included because they seem of independent interest.

\section{Non-squashing partitions }\label{Sec:NS}

Let $a(n)$ denote the number of non-squashing partitions of $n$.
It was shown by
Hirschhorn and Sellers\footnote{Hirschhorn and Sellers
regard the inequalities \eqref{EqNS} as purely arithmetic conditions
and do not mention stacking problems.}
\cite{HS03}
that $a(n)$ is equal to the number of ``binary partitions''
of $n$, that is, the number
of partitions of $n$ into powers of $2$.
See sequences
\htmladdnormallink{A000123}{http://www.research.att.com/cgi-bin/access.cgi/as/njas/sequences/eisA.cgi?Anum=A000123}
and
\htmladdnormallink{A018819}{http://www.research.att.com/cgi-bin/access.cgi/as/njas/sequences/eisA.cgi?Anum=A018819} in \cite{OEIS} 
for properties of the binary partition function
and references to the extensive literature.

In fact Hirschhorn and Sellers prove a more general result.
Let $s \ge 2$ be an integer.
Let us say that a partition \eqref{EqPart}
is {\em $s$-non-squashing} if
\beql{EqNS2}
(s-1)(p_1 + \cdots + p_j) \le p_{j+1} \mbox{~for~} 
1 \le j \le k-1 ~.
\eeq
If the $p_j$ are the labels of the boxes in a stack,
not only is no box squashed, no box
even comes within a factor of $s-1$ of being squashed.
A non-squashing partition as defined in the Introduction is
$2$-non-squashing.

\begin{thm}\label{thm1} (Hirschhorn and Sellers \cite{HS03}.)
The number
$a_s(n)$ of $s$-non-squashing partitions of $n$ is equal to
the number of ``$s$-ary'' partitions of $n$, that is,
the number of partitions of $n$ into
powers of $s$.
\end{thm}

\vspace{.4\baselineskip}

The following is an alternative proof of this result
which leads to a bijection between the two
families.

\begin{proof}
Let $a'_s(n)$ be the number of partitions of $n$
into powers of $s$, for some integer $s \ge 2$. Suppose
$$
n = s^{e_1} + s^{e_2} + \cdots + s^{e_l}
$$
is such a partition, where $n \ge s$.
If at least one of the parts is $1$ we can remove it and
obtain a partition of $n-1$ into powers of $s$;
if not, all the $e_i$ are
greater than $0$  and we can also divide by $s$ and obtain a partition of $n/s$.
Therefore $a'_s(n)$ satisfies the recurrence
\begin{eqnarray}\label{EqPS}
a'_s(n) & = & a'_s(n-1) \mbox{~if~} n \not\equiv 0 \mbox{~mod~} s ~, \nonumber \\
a'_s(n) & = & a'_s(n-1) + a'_s(n/s) \mbox{~if~} n \equiv 0 \mbox{~mod~} s ~,
\end{eqnarray}
for $n \ge s$.
The smallest $n$ for which there is a partition with more than
one part is $s$, so we have the initial conditions
\beql{EqPS2}
a'_s(0) = a'_s(1) = \cdots = a'_s(s-1) = 1 ~.
\eeq

On the other hand, let 
\beql{EqPart2}
n = p_1 + p_2 + \cdots + p_k \mbox{~with~}
1 \le p_1 \le p_2 \le \cdots \le p_k
\eeq
be an $s$-non-squashing partition of $n \ge s$.
If the largest part $p_k$ is strictly greater than
$\frac{(s-1)n}{s}$, then the sum of the other parts is 
strictly less than $\frac{n}{s}$, and we can subtract $1$ 
from the largest part and obtain 
an $s$-non-squashing partition of $n-1$.
(We omit the straightforward verification.)
If the largest part is equal to
$\frac{(s-1)n}{s}$ (implying $n \equiv 0 \mbox{~mod~} s$),
we can also delete the largest part and obtain 
an $s$-non-squashing partition of $n/s$.
Therefore $a_s(n)$ satisfies the recurrence
\begin{eqnarray}\label{EqNS3}
a_s(n) & = & a_s(n-1) \mbox{~if~} n \not\equiv 0 \mbox{~mod~} s ~, \nonumber \\
a_s(n) & = & a_s(n-1) + a_s(n/s) \mbox{~if~} n \equiv 0 \mbox{~mod~} s ~,
\end{eqnarray}
for $n \ge s$.
The smallest $n$ for which there is a partition with more than
one part is $s$ 
(where we have the partition with parts $1$ and $s-1$),
so we have the initial conditions
\beql{EqNS4}
a_s(0) = a_s(1) = \cdots = a_s(s-1) = 1 ~.
\eeq

Comparing \eqn{EqPS}, \eqn{EqPS2} with \eqn{EqNS3}, \eqn{EqNS4},
we conclude that $a'_{s}(n) = a_s(n)$
for all $n \ge 0$ and all $s \ge 2$,
which is the main result of \cite{HS03}.
\end{proof}

\vspace{.4\baselineskip}

\begin{figure}[htb]
   \begin{center}
      \epsfxsize=3.5in
      \leavevmode\epsffile{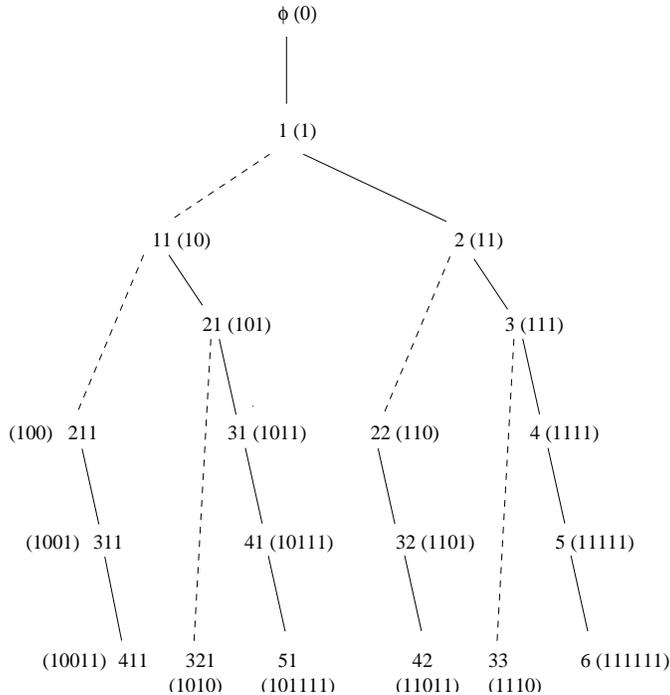}
    \end{center}
\caption{
Non-squashing partitions of the numbers $0, \ldots, 6$
arranged in tree structure.
The binary labels are shown in parentheses.
(Every node has out-degree $2$, but only edges between
partitions of $0, \ldots, 6$ are shown.) }
\label{Fig1}
\end{figure}

\begin{figure}[htb]
   \begin{center}
      \epsfxsize=3.5in
      \leavevmode\epsffile{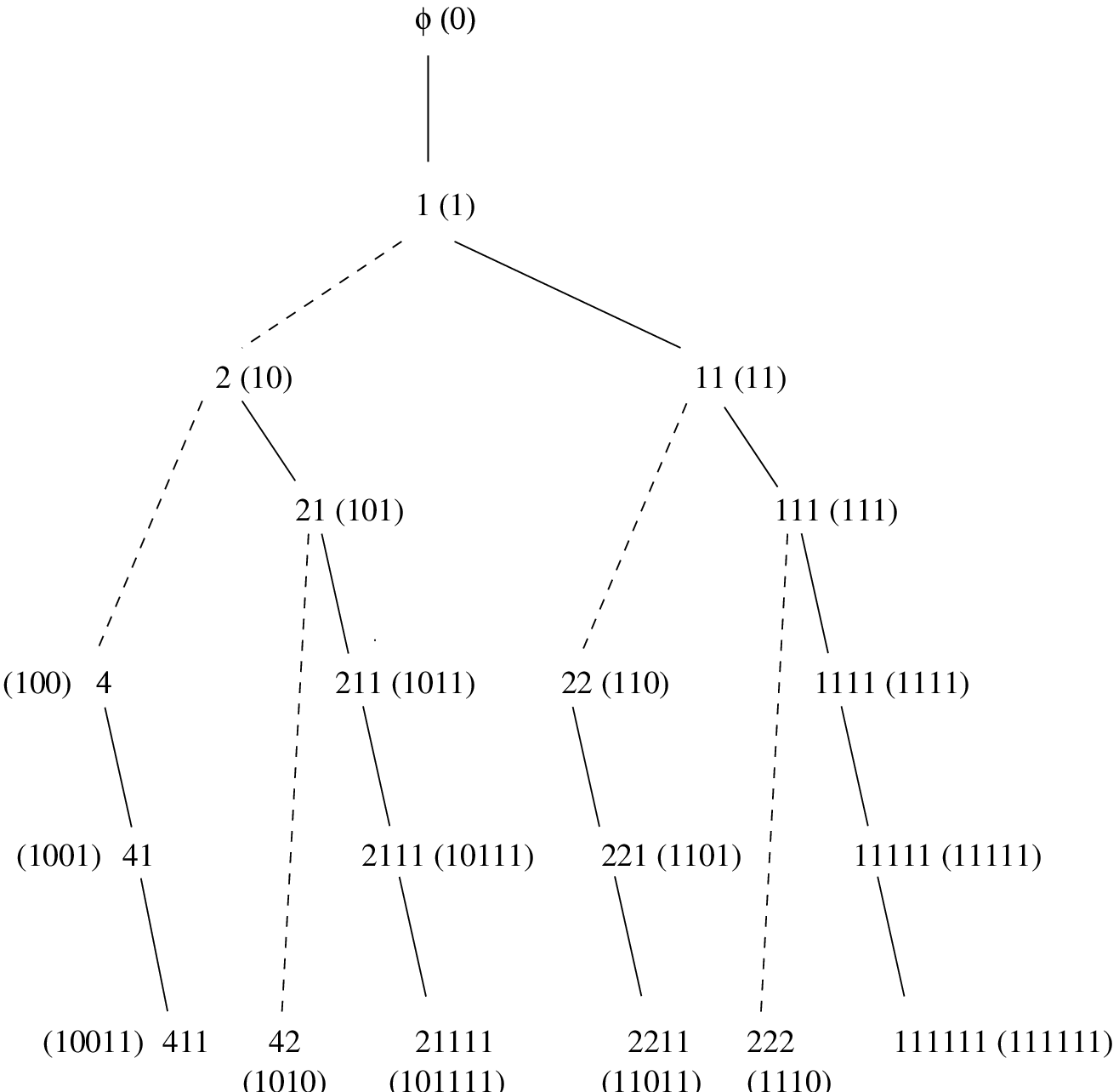}
    \end{center}
\caption{
Binary partitions of the numbers $0, \ldots, 6$
arranged in tree structure.
The binary labels are shown in parentheses.
(Every node has out-degree $2$, but only edges between
partitions of $0, \ldots, 6$ are shown.) }
\label{Fig2}
\end{figure}

The above proof associates each partition (from either family)
with a unique partition of a smaller number.
We can therefore arrange the partitions in each family into
a rooted tree, with the empty partition of $0$ as the root node.
Figures \ref{Fig1} and \ref{Fig2} show the beginnings of the two
trees for the case $s=2$.
(Most of the time we will adopt the standard convention
of writing partitions with the parts in nonincreasing order.)
Every node has two descendants and (except for the root)
one ancestor. We may label the edge leading from
a partition of $n/s$ to a partition of $n$ 
with $0$ (such edges are shown as broken lines
in Figs. \ref{Fig1} and \ref{Fig2}),
and the edge leading from
a partition of $n-1$ to a partition of $n$ 
with $1$ (the solid lines in the figures).

This associates a unique binary string with each partition in either tree.
A partition of $n$ in one tree receives the same
binary string as the corresponding partition of $n$ in the same
position in the other tree.
In this way we obtain a canonical numbering for the
$s$-non-squashing partitions,
a canonical numbering for the
partitions into powers of $s$,
and a bijection between them.

Table \ref{T1} shows the beginning of the bijection.
The first column gives the binary string $u$,
the second column gives the corresponding
$s$-non-squashing partition $P(u)$,
the third column gives the corresponding
$s$-ary partition $Q(u)$,
and the last column gives the number $n = n(u)$ that is 
partitioned by both $P(u)$ and $Q(u)$.

\begin{table}[htb]
\caption{Bijection between $s$-non-squashing partitions $P(u)$ and
$s$-ary partitions $Q(u)$; $u$ is the index (written as a binary number) and
$n$ is the number that is being partitioned.}
$$
\begin{array}{|r|l|l|r|} \hline
u & P(u) & Q(u) & n \\ \hline
0 & \emptyset & \emptyset & 0 \\
1 & 1 & 1 & 1 \\
10 & s-1,1 & s & s \\
11 & 2 & 1,1 & 2 \\
100 & s(s-1),s-1,1 & s^2 & s^2 \\ 
101 & s,1 & s,1 & s+1 \\ 
110 & 2(s-1),2 & s,s & 2s \\ 
111 & 3 & 1,1,1 & 3 \\ 
1000 & s^2(s-1),s(s-1),s-1,1 & s^3 & s^3 \\ 
1001 & s^2-s+1,s-1,1 & s^2,1 & s^2+1 \\ 
1010 & s^2-1,s,1 & s^2,s & s^2+s \\ 
1011 & s+1,1 & s,1,1 & s+2 \\ 
1100 & 2s(s-1),2(s-1),2 & s^2,s^2 & 2s^2 \\ 
1101 & 2s-1,2 & s,s,1 & 2s+1 \\ 
1110 & 3(s-1),3 & s,s,s & 3s \\ 
1111 & 4 & 1,1,1,1 & 4 \\ 
10000 & s^3(s-1),s^2(s-1),s(s-1),s-1,1 & s^4 & s^4 \\
10001 & s^3-s^2+1,s(s-1),s-1,1 & s^3,1 & s^3+1 \\
10010 & (s-1)(s^2+1),s^2-s+1,s-1,1 & s^3,s & s^3+s \\
10011 & s^2-s+2,s-1,1 & s^2,1,1 & s^2+2 \\
10100 & s(s^2-1),s^2-1,s,1 & s^3,s^2 & s^3+s^2 \\ 
10101 & s^2,s,1 & s^2,s,1 & s^2+s+1 \\
\hline
\end{array}
$$
\label{T1}
\end{table}

We note without proof the following properties of the bijection.

(i) For a nonzero string $u$, the number of parts in $P(u)$
is equal to $1$ plus the number of $0$'s in $u$, and
the number of parts in $Q(u)$
is equal to the number of $1$'s in $u$.

(ii) Thinking of $u$ now as the integer represented by the binary string,
the number $n = n(u)$ (given in the 
last column of the table) that is partitioned by both $P(u)$ and $Q(u)$
is defined by the recurrence
$$
n(0) = 0; ~ n(2u)=sn(u) \mbox{~for~} u \ge 1, ~ n(2u+1)=n(u) + 1 \mbox{~for~} u \ge 0
$$
(sequences 
\htmladdnormallink{A087808}{http://www.research.att.com/cgi-bin/access.cgi/as/njas/sequences/eisA.cgi?Anum=A087808},
\htmladdnormallink{A090639}{http://www.research.att.com/cgi-bin/access.cgi/as/njas/sequences/eisA.cgi?Anum=A090639},
etc. in  \cite{OEIS}). 

(iii) $P(u) = Q(u)$ if and only if $u = (4^k-1)/3$ for some
$k \ge 1$ (that is, if $u$ is the binary string $10101\ldots01$).

It is easy to go from the binary vector to the partitions and vice versa.
To obtain the $s$-non-squashing partition $P(u)$
corresponding to the binary vector $u$, we start with the empty
partition $P(u) = \emptyset$, and scan $u$ from left to
right (i.e. beginning with the most significant bit):
\begin{itemize}
\item
if we see a $1$, then if $P(u) = \emptyset$ set $P(u) = 1$, 
otherwise add $1$ to the largest part of $P(u)$,
\item
if we see a $0$, then if $P(u) = \emptyset$ set $P(u) = 0$, 
otherwise adjoin to $P(u)$ a part equal
to $s-1$ times the sum of the parts of $P(u)$.
\end{itemize}

Example for $s=3$. Suppose $u = 10110$. The successive terms in
the construction of $P(u)$ are
$$
\emptyset,~ 1,~ 2\,1,~ 3\,1,~ 4\,1,~ 10\,4\,1 ~.
$$

Likewise, to obtain the partition $Q(u)$ into powers of $s$,
again we start with the empty
partition $Q(u) = \emptyset$, and scan $u$ from left to right:
\begin{itemize}
\item
if we see a $1$, then append a part of size $1$ to $Q(u)$,
\item
if we see a $0$, then if $Q(u) = \emptyset$ do nothing,
otherwise multiply all the parts of $Q(u)$ by $s$.
\end{itemize}

Example for $s=3$. Again we take $u = 10110$.  The successive terms in
the construction of $Q(u)$ are
$$
\emptyset,~ 1,~ 3,~ 3\,1,~ 3\,1\,1,~ 9\,3\,3 ~.
$$
Thus the bijection associates these two
partitions of $15$: $P(u) = 10\,4\,1$ and $Q(u) = 9\,3\,3$.

Finally, we note that the numbers $a_s(n)$ have the generating function
\beql{GF1}
\sum_{n=0}^{\infty} a_s(n) \, x^n = \prod_{i=0}^{\infty} \frac{1}{1-x^{s^i}} ~.
\eeq

\section{Non-squashing partitions into distinct parts }\label{NSPD}
 From here on we consider only the case $s=2$, that is, 
non-squashing partitions.
One of the restrictions in the box-stacking
problem mentioned in the Introduction
is that the parts be distinct.
In this section we investigate the number $b(n)$
of non-squashing partitions of $n$ into distinct parts.
The first few values of $b(n)$ for $n = 0, 1, 2, \ldots$ are
\beql{EqNSD2aa}
1, \,1, \,1, \,2, \,2, \,3, \,4, \,5, \,6, \,7, \,9, \,10, \,13, \,14, \,18, \,19, \,24, \,25, \,31, \,32, \,40, \,41, \,50, \ldots
\eeq
(this is now sequence \htmladdnormallink{A088567}{http://www.research.att.com/cgi-bin/access.cgi/as/njas/sequences/eisA.cgi?Anum=A088567} in \cite{OEIS}).

\begin{thm}\label{thm2}
The numbers $b(n)$ satisfy the recurrence
\begin{eqnarray}\label{EqNSD2}
b(0) & = &  b(1) ~ = ~ 1 ~, \nonumber \\
b(2m) & = & b(2m-1) + b(m) - 1 \mbox{~for~} m \ge 1 ~, \nonumber \\
b(2m+1) & = & b(2m) + 1 \mbox{~for~} m \ge 1 ~. 
\end{eqnarray}
The generating function
$B(x) = \sum_{n=0}^{\infty} b(n) x^n$ satisfies
\beql{EqNSD2b}
B(x) = \frac{1}{1-x} B(x^2) - \frac{x^2}{1-x^2} ~,
\eeq
and is given explicitly by
\beql{EqNSD5}
B(x) ~=~
1 ~+~ \frac{x}{1-x} ~+~
\sum_{i=1}^{\infty}
\frac{x^{3\cdot 2^{i-1}}} {\prod_{j=0}^{i} (1-x^{2^j}) } ~.
\eeq
\end{thm}

\begin{proof}
We obtain the partitions
(non-squashing and with distinct parts being understood)
of an odd number $2m+1$ by adjoining a part of size $2m+1-j$
to a partition of $j$, for some $j = 0, 1, \ldots, m$
(since $2m+1-j > j$, these are indeed non-squashing).
Likewise we obtain the partitions
of $2m$ by adjoining a part of size $2m-j$
to a partition of $j$, for some $j = 0, 1, \ldots, m$,
except that if $j=m$ we cannot adjoin a part of size $m$ to
the partition consisting of a single $m$.
Thus we have
\begin{eqnarray}\label{EqNSD1}
b(0) & = &  b(1) ~ = ~ 1 ~, \nonumber \\
b(2m+1) & = & b(0) + b(1) + \cdots + b(m) \mbox{~for~} m \ge 1 ~, \nonumber \\
b(2m) & = & b(0) + b(1) + \cdots + b(m) - 1 \mbox{~for~} m \ge 1 ~,
\end{eqnarray}
from which \eqn{EqNSD2} follows.
After some algebra we find that \eqn{EqNSD2} implies
\beql{EqNSD2a}
B(x) = xB(x) + B(x^2) - \frac{x^2}{1+x} ~,
\eeq
and after rearranging we obtain \eqn{EqNSD2b}.
Equation \eqn{EqNSD2b} implies
$$
B(x^2) = \frac{1}{1-x^2} B(x^4) - \frac{x^4}{1-x^4} ~,
$$
and so on, and hence 
\begin{eqnarray}\label{EqNSD3}
B(x) &  ~=~ & \prod_{i=0}^{\infty} \frac{1}{1-x^{2^i}} ~-~
\sum_{i=1}^{\infty} \frac{x^{2^i}}
{\prod_{j=0}^{i-1} (1-x^{2^j}) \cdot ( 1+x^{2^{i-1}}) }  \nonumber \\
{ }     &  ~=~ & \prod_{i=0}^{\infty} \frac{1}{1-x^{2^i}} ~-~
\sum_{i=1}^{\infty} \frac{x^{2^i}}
{\prod_{j=0}^{i-2} (1-x^{2^j}) \cdot ( 1-x^{2^{i}}) } ~.
\end{eqnarray}
To simplify this we make use of an identity from \cite{HS03}:
if $m_1 < m_2 < \cdots < m_k$ are positive integers then
\beql{EqNSD4}
1 + \sum_{j=1}^{k}
\frac{x^{m_j}}
{(1-x^{m_1}) (1-x^{m_2}) \cdots (1-x^{m_j})}
~=~
\prod_{j=1}^{k} \frac{1}{1-x^{m_j}} ~.
\eeq
Applying this to the sum in \eqn{EqNSD3} and simplifying, we
eventually obtain \eqn{EqNSD5}.
\end{proof}

\begin{cor}
(i) The sequence $\{b(n)\}$ (see \eqn{EqNSD2aa}) has the property
that the sequence of partial sums 
\beql{EqNSD6}
1, \, 2, \, 3, \, 5, \, 7, \, 10, \, 14, \, 19, \, 25, \, 32, \, 41, \, 51, \, 64, \, 78, \, 96, \, 115, \, 139, \, 164, \ldots
\eeq
coincides with the odd-indexed subsequence $b(1), b(3), b(5), \ldots$.
The even-indexed subsequence $b(2), b(4), b(6), \ldots$
is obtained by adding $1$ to the terms of \eqn{EqNSD6}.
(ii) b(n), the number of
non-squashing partitions of n into distinct parts,
is equal to the number of partitions of
$n$ into powers of $2$ such that either all the parts
are equal to $1$ or, if the largest part
has size $2^i > 1$, then there is also at least one part
of size $2^{i-1}$.
\end{cor}

\begin{proof}
(i) The first assertion is equivalent to the algebraic identity
\beql{EqNSD7}
\frac{B(x)}{1-x} =
\frac{B(\sqrt{x}) - B(-\sqrt{x})}{2\sqrt{x}} ~,
\eeq
which is easily verified using 
\eqn{EqNSD2a} and \eqn{EqNSD3}.
The second assertion follows from \eqn{EqNSD2}.
Property (ii) is an immediate consequence of \eqn{EqNSD3}.
\end{proof}

Congruences satisfied by $a_s(n)$ have been
studied by many authors (see references in \cite{HS03}).
Here we record just one such result for $b(n)$.

\begin{cor}
The value of $b(n)$ mod $2$ is as follows (all congruences are mod $2$):
\beql{EqC1}
b(0) \equiv 1 ~,
\eeq
\beql{EqC2}
\mbox{if~}n \mbox{~is~odd}, ~   b(n) \equiv b(n-1) + 1 ~,
\eeq
\beql{EqC3}
b(8m+2) \equiv 1, ~  b(8m+6) \equiv 0 ~,
\eeq
\beql{EqC4}
b(16m+4) \equiv 0, ~  b(16m+12) \equiv 1  ~,
\eeq
\beql{EqC5}
\mbox{for~} m>0, ~ b(16m) \equiv b(8m), ~ b(32m+8) \equiv 0, ~ b(32m+24) \equiv 1 ~.
\eeq
For $m>0$, $b(8m)$ is the value of the bit immediately to the left
of the rightmost $1$ when $m$ is written in binary.
\end{cor}

\begin{proof}
\eqn{EqC2} follows from \eqn{EqNSD2}. To prove the first
assertion in \eqn{EqC3}, we repeatedly apply \eqn{EqNSD2}, obtaining
\begin{eqnarray}\label{EqC6}
b(8m+2)  & \equiv & b(8m+1) + b(4m+1) + 1 \nonumber \\
{ }      & \equiv & b(8m) + b(4m) + 1 \nonumber \\
{ }      & \equiv & b(8m-1) \nonumber \\
{ }      & \equiv & b(8m-2) + 1 \nonumber \\
{ }      & \equiv & \ldots \nonumber \\
{ }      & \equiv & b(8i-6) \nonumber \\
{ }      & \equiv & \ldots \nonumber \\
{ }      & \equiv & b(2) = 1  \nonumber
\end{eqnarray}
The other claims in \eqn{EqC3}--\eqn{EqC5}
are established in a similar way.
It is easily checked that the final assertion 
in the corollary is equivalent to \eqn{EqC5}.
(The final assertion
was discovered by noticing that the subsequence $\{b(8m)\}$
is, apart from the leading term, the same as
sequence \htmladdnormallink{A038189}{http://www.research.att.com/cgi-bin/access.cgi/as/njas/sequences/eisA.cgi?Anum=38189} in \cite{OEIS}.)
\end{proof}

\section{Non-squashing partitions by number of parts }\label{NSNOP}
Let $a(n,k)$ be the number of non-squashing partitions of $n$
into exactly $k$ parts. Table \ref{T2} shows the initial values
of this function.

\begin{table}[htb]
\caption{ Values of $a(n,k)$, the number of non-squashing partitions 
of $n$ into exactly $k$ parts. }
$$
\begin{array}{|r|rrrrr|} \hline
  & \multicolumn{5}{c|}{k} \\ 
n & 0 & 1 & 2 & 3 & 4 \\ \hline
0 & 1 & 0 & 0 & 0 & 0  \\
1 & 0 & 1 & 0 & 0 & 0  \\
2 & 0 & 1 & 1 & 0 & 0  \\
3 & 0 & 1 & 1 & 0 & 0  \\
4 & 0 & 1 & 2 & 1 & 0  \\ 
5 & 0 & 1 & 2 & 1 & 0  \\ 
6 & 0 & 1 & 3 & 2 & 0  \\ 
7 & 0 & 1 & 3 & 2 & 0  \\ 
8 & 0 & 1 & 4 & 4 & 1  \\ 
9 & 0 & 1 & 4 & 4 & 1  \\ 
10 & 0 & 1 & 5 & 6 & 2 \\ 
11 & 0 & 1 & 5 & 6 & 2 \\ 
12 & 0 & 1 & 6 & 9 & 4 \\ 
13 & 0 & 1 & 6 & 9 & 4 \\ 
\hline
\end{array}
$$
\label{T2}
\end{table}

\begin{thm}\label{thm3}
The numbers $a(n,k)$ satisfy the recurrence
\begin{eqnarray}\label{EqNOP1}
a(2m,k) & = & a(2m-1,k) + a(m,k-1) \mbox{~for~} m \ge 1, k \ge 1 ~, \nonumber \\
a(2m+1,k) & = & a(2m,k) \mbox{~for~} m \ge 1, k \ge 1 ~, 
\end{eqnarray}
with initial conditions
$$
a(0,0) = 1 \, ,  \,
a(n,0) = 0  \mbox{~for~} n \ge 1 \, , \,
a(n,k) = 0  \mbox{~for~} k > n \, , \,
a(n,1) = 1  \mbox{~for~} n \ge 1 ~ .
$$
In particular,
each odd-indexed row (except
for row $1$) in Table \ref{T2} is a copy of the previous row.
If the duplicate entries are omitted, the $k$-th column
has generating function
\beql{EqNOP2}
\sum_{k=0}^{\infty} a(2m,k) x^m ~=~
\frac{ x^{2^{k-2}} }
{(1-x) \cdot \prod_{j=0}^{k-2} (1-x^{2^j}) } ~,
\eeq
while if they are included we get the simpler expression
\beql{EqNOP3}
\sum_{k=0}^{\infty} a(m,k) x^m ~=~
\frac{ x^{2^{k-1}} }
{\prod_{j=0}^{k-1} (1-x^{2^j}) } ~.
\eeq
Equation \eqn{EqNOP3} implies that 
the number of non-squashing partitions of $n$ with $k$ parts
is equal 
(i) to the number of partitions of $n-2^{k-1}$ into powers 
of $2$ not exceeding $2^{k-1}$, and also (ii) to the
number of binary partitions of $n$ with largest part $2^{k-1}$.
\end{thm}

\begin{proof}
The recurrence \eqn{EqNOP1} follows at once from
the argument used to derive \eqn{EqNS3}.
The generating functions then follow from the recurrence; we omit the details.
\end{proof}

\vspace{.4\baselineskip}

For example, the $k=3$ column, omitting the odd-indexed terms, is
$$
0, \, 0, \, 1, \, 2, \, 4, \, 6, \, 9, \, 12, \, 16, \, 20, \, 25, \, 30, \, 36, \, 42, \, 49, \, 56, \, 64, \, 72, \, 81, \ldots ~,
$$
which is the sequence
of ``quarter-squares'', that is,
$$
a(2m,3) = \bigg\lfloor \frac{m}{2} \bigg\rfloor \bigg\lceil  \frac{m}{2} \bigg\rceil ~,
$$
with generating function
$$
\sum_{k=0}^{\infty} a(2m,3) x^m ~=~
\frac{x^2}{(1-x)^2(1-x^2)}
$$
(sequence \htmladdnormallink{A002620}{http://www.research.att.com/cgi-bin/access.cgi/as/njas/sequences/eisA.cgi?Anum=A002620}).

\section{Non-squashing partitions into $k$ distinct parts }\label{NSkDP}

Let $b(n,k)$ be the number of non-squashing partitions of $n$
into exactly $k$ distinct parts. Table \ref{T3} shows the initial values
of this function.

\begin{table}[htb]
\caption{ Values of $b(n,k)$, the number of non-squashing partitions 
of $n$ into exactly $k$ distinct parts. }
$$
\begin{array}{|r|rrrrr|} \hline
  & \multicolumn{5}{c|}{k} \\ 
n & 0 & 1 & 2 & 3 & 4 \\ \hline
0 & 1 & 0 & 0 & 0 & 0  \\
1 & 0 & 1 & 0 & 0 & 0  \\
2 & 0 & 1 & 0 & 0 & 0  \\
3 & 0 & 1 & 1 & 0 & 0  \\
4 & 0 & 1 & 1 & 0 & 0  \\ 
5 & 0 & 1 & 2 & 0 & 0  \\ 
6 & 0 & 1 & 2 & 1 & 0  \\ 
7 & 0 & 1 & 3 & 1 & 0  \\ 
8 & 0 & 1 & 3 & 2 & 0  \\ 
9 & 0 & 1 & 4 & 2 & 0  \\ 
10 & 0 & 1 & 4 & 4 & 0 \\ 
11 & 0 & 1 & 5 & 4 & 0 \\ 
12 & 0 & 1 & 5 & 6 & 1 \\ 
13 & 0 & 1 & 6 & 6 & 1 \\ 
14 & 0 & 1 & 6 & 9 & 2 \\ 
15 & 0 & 1 & 7 & 9 & 2 \\ 
16 & 0 & 1 & 7 & 12 & 4 \\ 
17 & 0 & 1 & 8 & 12 & 4 \\ 
\hline
\end{array}
$$
\label{T3}
\end{table}

Comparison of this table with Table \ref{T2} suggests that Table
\ref{T3} is obtained by displacing the $k$-th column
of Table \ref{T2} (for $k \ge 2$) downwards 
by $2^{k-2}$ positions. This is true, and we have:

\begin{thm}\label{thm4}
The numbers $b(n,k)$ satisfy 
\begin{eqnarray}\label{EqNSk1}
b(n,0) & = & a(n,0) \mbox{~for~} n \ge 0 ~, \nonumber \\
b(n,1) & = & a(n,1) \mbox{~for~} n \ge 0 ~, \nonumber \\
b(n,k) & = & a(n-2^{k-2},k) \mbox{~for~} n \ge 0, k \ge 2 ~.
\end{eqnarray}
Also
\beql{EqNSk7}
\sum_{k=0}^{\infty} b(n,k) x^n   =  
\frac{ x^{3 \cdot 2^{k-2}} }
{\prod_{j=0}^{k-1} (1-x^{2^j}) } \mbox{~for~} k \ge 2 ~.
\eeq
Equation \eqn{EqNSk7} implies that the number of non-squashing
partitions of $n$ with $k$ distinct parts is equal to the number
of partitions of $n-3 \cdot 2^{k-2}$ into powers of $2$ not
exceeding $2^{k-1}$.
\end{thm}

\begin{proof}
For this discussion we write the parts in nondecreasing order.
The non-squashing partition (of some very large number)
having the slowest growth 
begins
\beql{EqNSk4}
1, \,  1, \,  2, \,  4, \,  8, \,  16, \,  32, \,  64, \ldots ~, \, 
\eeq
while the non-squashing partition with distinct parts and the slowest growth 
is the sequence $\gamma(1), \gamma(2), \gamma(3),  \ldots $
given by
\beql{EqNSk5}
1, \,  2, \,  3, \,  6, \,  12, \,  24, \,  48, \,  96, \ldots 
\eeq
with $\gamma(i)=i$ for $i \le 3$,
$\gamma(i)=3\cdot 2^{i-3}$ for $i \ge 3$.
The difference between \eqn{EqNSk4} and \eqn{EqNSk5} is
\beql{EqNSk6}
0, \,  1, \,  1, \,  2, \,  4, \,  8, \,  16, \,  32, \ldots ~.
\eeq
One can now verify that adding the initial $k$ terms of \eqn{EqNSk6}
term-by-term to the parts of a non-squashing partition of
$n$ into $k$ parts provides a bijection with a
non-squashing partition of $n$ into $k$ distinct parts,
and establishes the relations in \eqn{EqNSk1}.

For example,
the non-squashing partitions of $n=4, \ldots,8$ into $k=3$ parts are:
\begin{eqnarray}\label{EqNSk2}
4 & : & 112  \nonumber \\
5 & : & 113 \nonumber \\
6 & : & 114, 123 \nonumber \\
7 & : & 115, 124 \nonumber \\
8 & : & 116, 125, 134, 224 ~.
\end{eqnarray}
On the other hand,
the non-squashing partitions of $n=6, \ldots,10$ into $k=3$ distinct parts are:
\begin{eqnarray}\label{EqNSk3}
6 & : & 123  \nonumber \\
7 & : & 124 \nonumber \\
8 & : & 125, 134 \nonumber \\
9 & : & 126, 135 \nonumber \\
10 & : & 127, 136, 145, 235 ~.
\end{eqnarray}
\noindent
Adding $0, 1, 1$ term-by-term to the partitions in \eqn{EqNSk2}
yields the partitions in \eqn{EqNSk3}.

The generating function \eqn{EqNSk7} now follows from
\eqn{EqNOP3} and \eqn{EqNSk1}.

\end{proof}

\section{Non-squashing partitions into distinct parts with largest part $m$ }\label{NSLP}
Let $c(n,k)$ be the number of non-squashing partitions of $n$
into distinct parts of which the greatest is $m$.
Table \ref{T4} shows the initial values.

\begin{table}[htb]
\caption{ Values of $c(n,k)$, the number of non-squashing partitions 
of $n$ into distinct parts of which the greatest is $m$
(the blank entries are zero).}
$$
\begin{array}{|r|rrrrrrrrrrrrr|} \hline
  & \multicolumn{13}{c|}{m} \\ 
n & 0 & 1 & 2 & 3 & 4 & 5 & 6 & 7 & 8 & 9 & 10 & 11 & 12 \\ \hline
0 & 1 &   &   &   &   &   &   &   &   &   &    &    &    \\
1 & 0 & 1 &   &   &   &   &   &   &   &   &    &    &    \\
2 & 0 & 0 & 1 &   &   &   &   &   &   &   &    &    &    \\
3 & 0 & 0 & 1 & 1 &   &   &   &   &   &   &    &    &    \\
4 & 0 & 0 & 0 & 1 & 1 &   &   &   &   &   &    &    &    \\
5 & 0 & 0 & 0 & 1 & 1 & 1 &   &   &   &   &    &    &    \\
6 & 0 & 0 & 0 & 1 & 1 & 1 & 1 &   &   &   &    &    &    \\
7 & 0 & 0 & 0 & 0 & 2 & 1 & 1 & 1 &   &   &    &    &    \\
8 & 0 & 0 & 0 & 0 & 1 & 2 & 1 & 1 & 1 &   &    &    &    \\
9 & 0 & 0 & 0 & 0 & 0 & 2 & 2 & 1 & 1 & 1 &    &    &    \\
10 & 0 & 0 & 0 & 0 & 0 & 2 & 2 & 2 & 1 & 1 & 1  &    &    \\
11 & 0 & 0 & 0 & 0 & 0 & 0 & 3 & 2 & 2 & 1 & 1  & 1  &    \\
12 & 0 & 0 & 0 & 0 & 0 & 0 & 3 & 3 & 2 & 2 & 1  & 1  & 1  \\
\hline
\end{array}
$$
\label{T4}
\end{table}

\begin{thm}\label{thm5}
(i) The nonzero values of $c(n,m)$ lie within a certain strip:
$$
c(n,m) = 0 \mbox{~if~} m < n/2 \mbox{~or~if~} n < m ~.
$$
(ii) For $m \le n \le 2m$,
\beql{EqLP1}
c(n,m) = \sum_{i=0}^{m-1} c(n-m, i) ~.
\eeq
(iii) For $m \le n \le 2m$,
\begin{eqnarray}\label{EqLP2}
c(n,m)  & = & b(n-m) \mbox{~if~} n < 2m ~, \nonumber \\
        & = & b(n-m) - 1 \mbox{~if~} n = 2m ~.
\end{eqnarray}
\end{thm}

\begin{proof}
(i) The slowest-growing non-squashing partition into distinct parts is
\eqn{EqNSk5}, so no partition can have $n > 2m$. The second
assertion is immediate from the definition of $c(n,m)$.

(ii) This is a consequence of the fact that removing the largest
part leaves a partition with largest part $ \le m-1$.

(iii) When the largest part is removed,
we obtain a non-squashing partition of $n-m$ into distinct parts.
Conversely, given a 
non-squashing partition of $n-m$ into distinct parts,
we obtain a partition of $n$ with largest part $m$
by adjoining a part of size $m$, with the single exception
that we cannot adjoin a part of size $m$ to
the partition consisting of a single
part of size $m$.

\end{proof}

\section{Solution to the box-stacking problem }\label{BOXES}
We can now give the solution to the box-stacking
problem mentioned in the Introduction.

\begin{thm}\label{thm7a}
There is a bijection between
non-squashing stacks of boxes in which the largest box has
label $n$ and non-squashing partitions of $2n$ into
distinct parts, i.e.
\beql{Eq12}
f(n) - f(n-1) = b(2n) ~.
\eeq
\end{thm}

\begin{proof}
Let
$$
1 \le p_1 < p_2 < \ldots < p_k = n
$$
be a non-squashing stack of boxes in which the largest box has
label $n$.
Let $r = p_1 + \cdots + p_{k-1}$ (take $r=0$ if $k=1$). Then 
$r \le p_k = n$.
If we increase the largest part by $n-r$ we obtain
a non-squashing partition of $2n$.
Conversely, suppose $1 \le p_1 < p_2 < \ldots < p_k$
is a non-squashing partition of $2n$ into distinct parts.
Let $r = p_1 + \cdots + p_{k-1}$. Then 
$r+p_k = 2n$, $r <p_k$, which implies $r < n < p_k$.
So we may reduce the largest part to $n$,
obtaining a non-squashing stack with largest part labeled $n$.
\end{proof}

Equation \eqn{Eq12}
could also be derived from the fact that
$$
f(n) = \sum_{i=0}^{\binom{n+1}{2}} \sum_{j=0}^{n} c(i,j)
     = \sum_{i=0}^{2n} \sum_{j=0}^{n} c(i,j) ~.
$$
\begin{cor}\label{thm7}
The numbers $f(n)$ have generating function
\beql{EqB1}
F(x) = \sum_{n=0}^{\infty} f(n) x^n = \frac{B(x)-x}{(1-x)^2} ~,
\eeq
where $B(x)$ is given in Theorem \ref{thm2}.
Also, $F(x)$ satisfies
\beql{EqB1a}
F(x) ~=~
\frac{ (1+x)^2}{1-x} F(x^2) ~-~
\frac{x(1-2x^2)}{(1-x)^2(1-x^2)} ~.
\eeq
\end{cor}

\vspace{.4\baselineskip}

\begin{proof}
 From Theorem \ref{thm7a} we know that
$$
f(n) = b(0) + b(2) + \cdots + b(2n) ~,
$$
so
$$
F(x) ~=~
\frac{1}{1-x} \, 
\frac{ B(\sqrt{x}) + B(-\sqrt{x})}{2} ~.
$$
So \eqn{EqB1} will follow if we can show that
$$
\frac{2(B(x)-x)}{1-x} ~=~ B(\sqrt{x}) + B(-\sqrt{x}) ~.
$$
However, from  \eqn{EqNSD7} we know that
$$
\frac{2 \sqrt{x} B(x)}{1-x} ~=~  B(\sqrt{x}) - B(-\sqrt{x}) ~.
$$
So we must show that
$$
B(\sqrt{x}) ~=~ \frac{B(x)-x}{1-x} ~+~ \frac{\sqrt{x} B(x)}{1-x} ~,
$$
which follows immediately from \eqn{EqNSD2b}.
Equation \eqn{EqB1a} then follows using \eqn{EqNSD2a}.
\end{proof}

\vspace{.4\baselineskip}

The first few values of $f(n)$ for $n = 0, 1, 2, \ldots$ are
\beql{EqB2}
1, \, 2, \, 4, \, 8, \, 14, \, 23, \, 36, \, 54, \, 78, \, 109, \, 149, \, 199, \, 262, \, 339, \, 434, \, 548, \, 686,  \ldots
\eeq
(sequence \htmladdnormallink{A089054}{http://www.research.att.com/cgi-bin/access.cgi/as/njas/sequences/eisA.cgi?Anum=A089054}).

\vspace{.4\baselineskip}

The original version of the problem had $n+1$ boxes
labeled $0, 1, \ldots, n$.  Since the box labeled $0$ 
may be included in the stack or not, without 
changing the non-squashing property, the answer
to this problem is $2f(n)$.

\section{Stacks with a given number of boxes }\label{BOXES2}
In this final section we determine the numbers $f(n,k)$,
the number of non-squashing stacks of boxes in
which the largest box has label $\le n$
and there are exactly $k$ boxes in the stack.
Table \ref{T5} shows the initial values
of this function.

\begin{table}[htb]
\caption{ Values of $f(n,k)$, the number of stacks
in which there are exactly $k$ boxes
and the largest box is $\le n$. }
$$
\begin{array}{|r|rrrrrr|} \hline
   &  \multicolumn{6}{c|}{k} \\
n & 0 & 1 & 2 & 3 & 4 & 5 \\ \hline
0 & 1 & 0 & 0 & 0 & 0 & 0 \\
1 & 1 & 1 & 0 & 0 & 0 & 0 \\
2 & 1 & 2 & 1 & 0 & 0 & 0\\
3 & 1 & 3 & 3 & 1 & 0 & 0\\
4 & 1 & 4 & 6 & 3 & 0 & 0\\
5 & 1 & 5 & 10 & 7 & 0 & 0 \\
6 & 1 & 6 & 15 & 13 & 1 & 0 \\
7 & 1 & 7 & 21 & 22 & 3 & 0 \\
8 & 1 & 8 & 28 & 34 & 7 & 0 \\
9 & 1 & 9 & 36 & 50 & 13 & 0\\
10 & 1 & 10 & 45 & 70 & 23 & 0 \\
11 & 1 & 11 & 55 & 95 & 37 & 0 \\
12 & 1 & 12 & 66 & 125 & 57 & 1 \\
\hline
\end{array}
$$
\label{T5}
\end{table}

\begin{thm}\label{thm8}
We have $f(n,0) = 1$ for all $n$,
and for $n \ge 1$, $k \ge 1$,
\beql{EqB8}
f(n,k) = \sum_{p=0}^{\min\{k-1, n-\gamma(k)\}}
~ \sum_{m=p}^{n-\gamma(k)} \,
(n-\gamma(k)+1-m) \, a(m,p) ~.
\eeq
\end{thm}

\begin{proof}
We first determine $f(n,k)-f(n-1,k)$, that is,
the number of stacks
$$
1 \le p_1 < p_2 < \ldots < p_k = n
$$
in which the largest box is labeled $n$.
Let $q_i = p_i - \gamma(i)$  for $i=1, \ldots, k$
(cf. \eqn{EqNSk5}),
so that
$$
0 \le q_1 \le q_2 \le \ldots \le q_k = n-\gamma(k) ~.
$$
Some of the $q_i$ may be zero. The nonzero elements 
among $q_1, \ldots, q_{k-1}$ (if any)
form a non-squashing partition into $p$ parts of
some number $m$ between
$0$ and $q_k$, where $0 \le p \le k-1$.  Hence
\beql{Eq831}
f(n,k)-f(n-1,k)
~=~
\sum_{m \le n - \gamma(k) } ~ \sum_{p \le k-1} \, a(m,p) ~,
\eeq
and so
\beql{Eq832}
f(n,k)
~=~
\sum_{p=0}^{k-1} ~ \sum_{\tau = k}^{n} ~ \sum_{m = p}^{\tau - \gamma(k)} \, a(m,p) ~.
\eeq
Equation \eqn{EqB8} follows when we collect terms.
\end{proof}
\vspace{.4\baselineskip}

\end{document}